\documentclass[11pt]{amsart}

\usepackage{enumerate}
\usepackage{latexsym}
\usepackage[pdftex]{graphicx}
\usepackage{amssymb}
\usepackage{xcolor}
\usepackage{stackrel}
\usepackage[normalem]{ulem}

\newtheorem{theorem}{Theorem}
\newtheorem{lemma}{Lemma}
\newtheorem{proposition}{Proposition}
\newtheorem{remark}{Remark}
\newtheorem{example}{Example}

\newtheorem{definition}{Definition}
\newtheorem{corollary}{Corollary}

\newtheorem{problem}{Problem}
\newtheorem*{problem*}{Problem}

\newcommand{\beq}{\begin{equation}}
\newcommand{\eeq}{\end{equation}}
\newcommand{\beqna}{\begin{eqnarray*}}
\newcommand{\eeqna}{\end{eqnarray*}}
\newcommand{\beqn}{\begin{equation*}}
\newcommand{\eeqn}{\end{equation*}}
\newcommand{\bp}{\begin{proof}}
\newcommand{\ep}{\end{proof}}
\newcommand{\bprop}{\begin{proposition}}
\newcommand{\eprop}{\end{proposition}}
\newcommand{\bt}{\begin{theorem}}
\newcommand{\et}{\end{theorem}}
\newcommand{\bex}{\begin{example}}
\newcommand{\eex}{\end{example}}
\newcommand{\bc}{\begin{corollary}}
\newcommand{\ec}{\end{corollary}}
\newcommand{\bl}{\begin{lemma}}
\newcommand{\el}{\end{lemma}}
\newcommand{\bprob}{\begin{problem}}
\newcommand{\eprob}{\end{problem}}
\newcommand{\br}{\begin{remark}}
\newcommand{\er}{\end{remark}}
\newcommand{\bd}{\begin{definition}}
\newcommand{\ed}{\end{definition}}

\begin{document}

\title
[On  bodies floating in equilibrium in every orientation ]
{On  bodies floating in equilibrium in every orientation}

\author[D. Ryabogin]{Dmitry Ryabogin}
\address{Department of Mathematics, Kent State University,
Kent, OH 44242, USA} \email{ryabogin@math.kent.edu}

\thanks{The   author is  supported in
	part by   Simons Collaboration Grant for Mathematicians program 638576,
	by U.S.~National Science Foundation Grant
	DMS-1600753 and by United States - Israel Binational Science Foundation (BSF)}

\keywords{Floating bodies, Ulam's problem, normal curvature}

\begin{abstract}
   Ulam's problem  19 from the Scottish Book asks:   {\it  is a solid of uniform density which floats in water in every position  necessarily a sphere?}  We  obtain several results related to this problem.

\end{abstract}

\maketitle

\section{Introduction}

Let the density of water be $1$ and assume that a convex body $K\subset {\mathbb R^3}$ of uniform density ${\mathcal D}\in(0,1)$ is submerged into water. We say that $K$ floats in equilibrium in the direction $\xi$ orthogonal to the water surface if the line  $\ell(\xi)$  connecting the center of mass of $K$ and the center of mass of the submerged part is parallel to $\xi$. We say that
$K$ floats in equilibrium in every orientation  if $\ell(\xi)$ is parallel to $\xi$ for every $\xi$.

The following intriguing problem was  proposed by Ulam \cite[Problem 19]{U}: {\it If a convex body  $K\subset {\mathbb R^3}$ made of material of uniform density ${\mathcal D}\in(0,1)$  floats in equilibrium in any orientation in water, must $K$ be spherical? }

Schneider \cite{Sch1} and Falconer   \cite{Fa}  showed that this is true, provided $K$ is centrally symmetric and ${\mathcal D}=\frac{1}{2}$. However, it has been recently proven in \cite{R2} that there are non-centrally-symmetric convex bodies of density ${\mathcal D}=\frac{1}{2}$ that float in equilibrium in every orientation.

The ``two-dimensional version" of the problem is also very interesting. In this case, we consider  floating logs of uniform cross-section, and seek for the ones that will float in every orientation with the axis horizontal. In other words, our cross-section  $K$ is a convex set in ${\mathbb R^2}$ and the water surface is a line that cuts off a set of the given area from $K$.
If ${\mathcal D}=\frac{1}{2}$, Auerbach  \cite{A} has exhibited logs with non-circular cross-section, both convex and non-convex, whose  boundaries  are   so-called Zindler curves \cite{Zi}.
More recently,    Bracho,  Montejano and   Oliveros  \cite{BMO} showed that 
for densities ${\mathcal D}=\frac{1}{3}$, $\frac{1}{4}$, $\frac{1}{5}$ and  $\frac{2}{5}$ the answer is affirmative, while  Wegner 
proved that for some  other  values of ${\mathcal D}\neq \frac{1}{2}$ the answer is negative, \cite{Weg1}, \cite{Weg2}; see also related  results of V\'arkonyi \cite{V1}, \cite{V2}. Overall, the case of general ${\mathcal D}\in (0,1)$ is notably involved and widely open.

No results in ${\mathbb R^3}$ are known for  densities ${\mathcal D}\in (0,1)$ different from $\frac{1}{2}$ and no counterexamples have been found so far. 
In this paper we prove and recall several results which  were used in the case of density $\frac{1}{2}$, \cite{R2}, and which, we believe,  would help to attack the problem for other densities.
We begin with
\bt\label{Dpr}
Let $d\ge 3$, let $K\subset {\mathbb R^d}$ be a convex body  and let $\delta\in (0,\textnormal{vol}_d(K))$.

If $K$ floats in equilibrium at the level $\delta$ in every orientation, then, for all  hyperplanes $H$ that cut off  the parts of volume $\delta$ from $K$,  the cutting sections $K\cap H$ have equal moments of inertia with respect to all $(d-2)$-dimensional planes $\Pi\subset H$ passing through the center of mass of $K\cap H$ and these moments are independent of $H$ and $\Pi$.

Conversely, let $K$ have a $C^1$-smooth boundary and let the center of mass of $K$ coincide with the center of mass of the surface of centers, i.e.,  the locus of the centers of mass of  all parts of volume $\delta$ that are cut off by the  cutting hyperplanes $H$. If  all cutting sections $K\cap H$ have equal moments of inertia with respect to all $(d-2)$-dimensional planes $\Pi\subset H$ passing through the center of mass of $K\cap H$ and these moments are independent of $H$ and $\Pi$, then $K$ floats in equilibrium at the level $\delta$ in every orientation.
\et
This  Theorem \footnote{This  result  was also  recently obtained  in \cite[Theorem 1.1]{FSWZ}, but  the case $\delta= \frac{\textrm{vol}_d(K)}{2}$
is considered under the  assumption that the Dupin floating body coincides with the B\'ar\'any-Larman-Sh\"utt-Werner floating body and it 
 is  a single point.}  gives an affirmative  answer to a question mentioned in  \cite[page 20, line 14 from below]{CFG}: ``It seems that the floating body problem is just (V, I)".
An analogous    Theorem  for  $d=2$ was obtained by Davidov   \cite{Da} and independently by Auerbach  \cite{A}, see Theorem 
\ref{Dakr} and Remark  \ref{nun} at the end of  Section \ref{Au}.

\bc\label{equi}
Let $d\ge 3$, 
let a convex body $K$ have a $C^1$-smooth boundary  and let 
$\delta\in (0, \textnormal{vol}_d(K))$. 
Assume also that  the center of mass of $K$ coincides with the center of mass of the surface of centers. 
If for  every  hyperplane $H$ that cuts off  the part of volume $\delta$ from $K$
every cutting section $K\cap H$  is $(d+1)$-equichordal, i.e., if there exists a constant $c$ such that for every line $l\subset 
K\cap H$ passing through the center of mass ${\mathcal C}(K\cap H)$  and having two points of intersection $\zeta_{\pm}(l)$ with the boundary of $K$ one has 
$$
\textnormal{dist}^{d+1}({\mathcal C}(K\cap H),\zeta_{+}(l))+\textnormal{dist}^{d+1}({\mathcal C}(K\cap H),\zeta_{-}(l))=c,
$$
 then $K$ floats in equilibrium in every orientation. 
\ec

Using the results in \cite{R1} and \cite{R2} one can show that
 the converse is not true, provided $\delta=\frac{\textnormal{vol}_d(K)}{2}$, i.e., there exists a non-centrally-symmetric body of revolution $K$ that  floats in equilibrium in every orientation, yet not every section $K\cap H$ by the hyperplane that cuts off the part of volume $\delta$ is $(d+1)$-equichordal.  On the other hand, it was proved in  \cite{R1}  that if $K$ is a body of revolution, then the condition 
that $K\cap H$ is $(d+1)$-equichordal for every hyperplane $H$  that cuts off  the part of volume $\delta$ from 
$K$ yields that it is the Euclidean ball.

\bprob\label{abdul1}
Is it possible to construct  a convex body $K$ and find $\delta\in (0, \textrm{vol}_d(K))$, $\delta\neq \frac{\textnormal{vol}_d(K)}{2}$,
so that $K\cap H$  is $(d+1)$-equichordal for every  hyperplane $H$ that cuts off  the part of volume $\delta$ from $K$, but 
$K$ is not an Euclidean ball?   
\eprob

We refer the reader to  \cite[pgs. 9-11]{CFG},  \cite[Chapter 6]{Ga}  and references therein for the information about equichordal bodies.

We also have

\bc\label{Al} 
Let $d\ge 2$ and  let a sequence $(\delta_n)_{n=1}^{\infty}$ of positive numbers  be such that the Dupin floating body $K_{[{\delta_n}]}$ coincides with the floating body $K_{\delta_n}$ for all $n\in {\mathbb N}$ and $\delta_n\to 0$ as $n\to\infty$.
If  $K$ floats in equilibrium in every orientation for all  levels $\delta_n$,   then $K$ is a Euclidean ball.      
\ec

Using Theorem \ref{Dpr} and the results of Myroshnychenko and Saroglou \cite{MRS},  one can also give a different  proof \footnote{See \cite[Theorem 1.2]{FSWZ} for a third  proof of this statement.}  of the  aforementioned result of Schneider and Falconer  obtained in \cite{Sch1} and \cite{Fa}  via spherical harmonics.

\bt\label{CS}
Let $d\ge 3$ and let $K\subset {\mathbb R^d}$ be a centrally-symmetric convex body. 
If $K$ floats in equilibrium in every orientation at the level $\delta=\frac{\textnormal{vol}_d(K)}{2}$ 
then 
$K$ is a Euclidean ball.
\et

Most of the results of this paper, as well as many  other results on floating bodies,  follow from the classical theorems of Dupin which, we believe, were missed by the mathematical community, \cite[Chapter XXIV]{DVP}, \cite[Hydrostatics, Part I]{Zh}).  In Sections \ref{El} and \ref{DP} we formulate  and  prove  these theorems in ${\mathbb R^d}$, $d\ge 3$ (see also \cite[Appendices A and B]{R2}).  

We refer the interested  reader to  \cite[pgs. 90-93]{M}, \cite[pgs. 19-20]{CFG},  \cite[pgs. 376-377]{Ga},  \cite[pgs. 560-563]{Sch2}, and \cite{G},  for an exposition of known results related to Ulam's Problem 19;   see also \cite{O}, \cite{Od}, \cite{HSW}, \cite{KO}, \cite{Gr} and  \cite{Mo} for related results.
The  {\it floating body problems}  appear in several areas of mathematics and, among other things, are related to the  Busemann-Petty problems in asymptotic geometric analysis \cite{BP},  to problems in statistics \cite{NSW}, and to
problems about polytopal approximation, \cite{B}, \cite{BL}, \cite{S2},  
\cite{BLW}.
We also refer the reader to \cite{MR}, \cite{St}, \cite{S1},  \cite{SW1},  \cite{SW2}, \cite{W},  and references therein for other works on floating bodies.

The paper is structured as follows. In the next section we recall some well-known facts about floating bodies and formulate the Theorems of Dupin 
in ${\mathbb R^d}$, $d\ge 3$.  We  prove  these theorems  in Section \ref{DP}. The proofs of   Lemma \ref{tr}, Theorems \ref{Dpr} and  \ref{CS}, and Corollaries  \ref{equi} and \ref{Al}  are given in Section \ref{Au}.

\section{Notation, basic definitions and Theorems of Dupin}\label{El}

\subsection{Notation and basic definitions}
 A convex body $K\subset {\mathbb R^d}$, $d\ge 2$,   is a convex compact set with a non-empty interior $\textrm{int} K$. 
 The boundary of $K$ is denoted by
 $\partial K$. We say that $K$ is  strictly convex  
 if $\partial K$ does not contain a segment.
 We say that $K$ is origin-symmetric if $K=-K$ and centrally-symmetric if there exists $p\in{\mathbb R^d}$ such that $K-p=\{q-p:\,q\in K\}$ is origin-symmetric.
 For $d\ge 2$ we denote by $S^{d-1}$  the unit sphere in ${\mathbb R^d}$ centered at the origin.  Given  $\xi\in S^{d-1}$ we denote by
  $\xi^{\perp}=\{p\in{\mathbb R^d}:\, p\cdot \xi=0\}$ the    subspace orthogonal to $\xi$, where   $p\cdot\xi=p_1\xi_1+\dots +p_d\xi_d$ is a usual inner product in ${\mathbb R^d}$.  
 The symbol  $``+"$ stands for the usual Minkowski (vector) addition, i.e., given two sets $D$ and $E$ in ${\mathbb R^d}$, $D+E=\{d+e:\,d\in D,\,e\in E   \}$. 
 Let $W_j$ be a $j$-dimensional  plane in ${\mathbb R^d}$, $1\le j\le d$.
 The {\it center of mass} of a  compact convex set $K\subset W_j$ with a non-empty relative interior will be denoted by ${\mathcal C}(K)$,
 $$
 {\mathcal C}(K)=\frac{1}{\textrm{vol}_j(K)}\int\limits_{K}xdx,
 $$
 where $\textrm{vol}_j(K)$ is the $j$-dimensional volume of $K$ in ${\mathbb R^j}$.
 We say that a  hyperplane $H$ is the supporting hyperplane of  a convex body $K$ if $K\cap H\neq \emptyset$, but $\textrm{int}\,K\cap H=\emptyset$. 
 
  If   $K\subset {\mathbb R^d}$ is a convex body  containing a point $p$  in its   interior, the {\it radial function} of $K$ with respect to $p$  in the direction $\theta\in S^{d-1}$   is defined as 
  $$
  \rho_{K,\,p}(\theta)=\max\{\lambda>0:\,p+\lambda \theta\in K\}.
  $$
  In particular, if $p$ is the origin, we will use the notation $$\rho_{K}(\theta)=\max\{\lambda>0:\lambda \theta\in K\}.
  $$
 Let  $m\in{\mathbb N}$. We say that  a convex body $K$ is of class $C^m({\mathbb R^d})$ (or $K$ has a $C^m$-smooth boundary) if for every point $z$ on the boundary $\partial K$ of $K\subset {\mathbb R^d}$ there exists a neighborhood $U_z$ of $z$  in ${\mathbb R^d}$ such that  $\partial K\cap U_z$ can be written as a graph of a function having 
all continuous partial derivatives  up to  the $m$-th order. 
The {\it Hausdorff distance} between two convex bodies $K$ and $L$ is defined as 
$$
d(K,L)=\sup\limits_{\{\theta\in S^{d-1}\}}|h_K(\theta)-h_L(\theta)|,
$$
where
$h_K$, $h_L$ are the {\it support functions} of bodies $K$, $L$, and for any $\theta\in S^{d-1}$,
$h_K(\theta)=\sup\limits_{\{y\in K \}}\theta\cdot y$.
A symbol $\,\square\,$ denotes  end of the proof.

\medskip

\begin{figure}[h]
	\centering
	\includegraphics[height=2.8in]{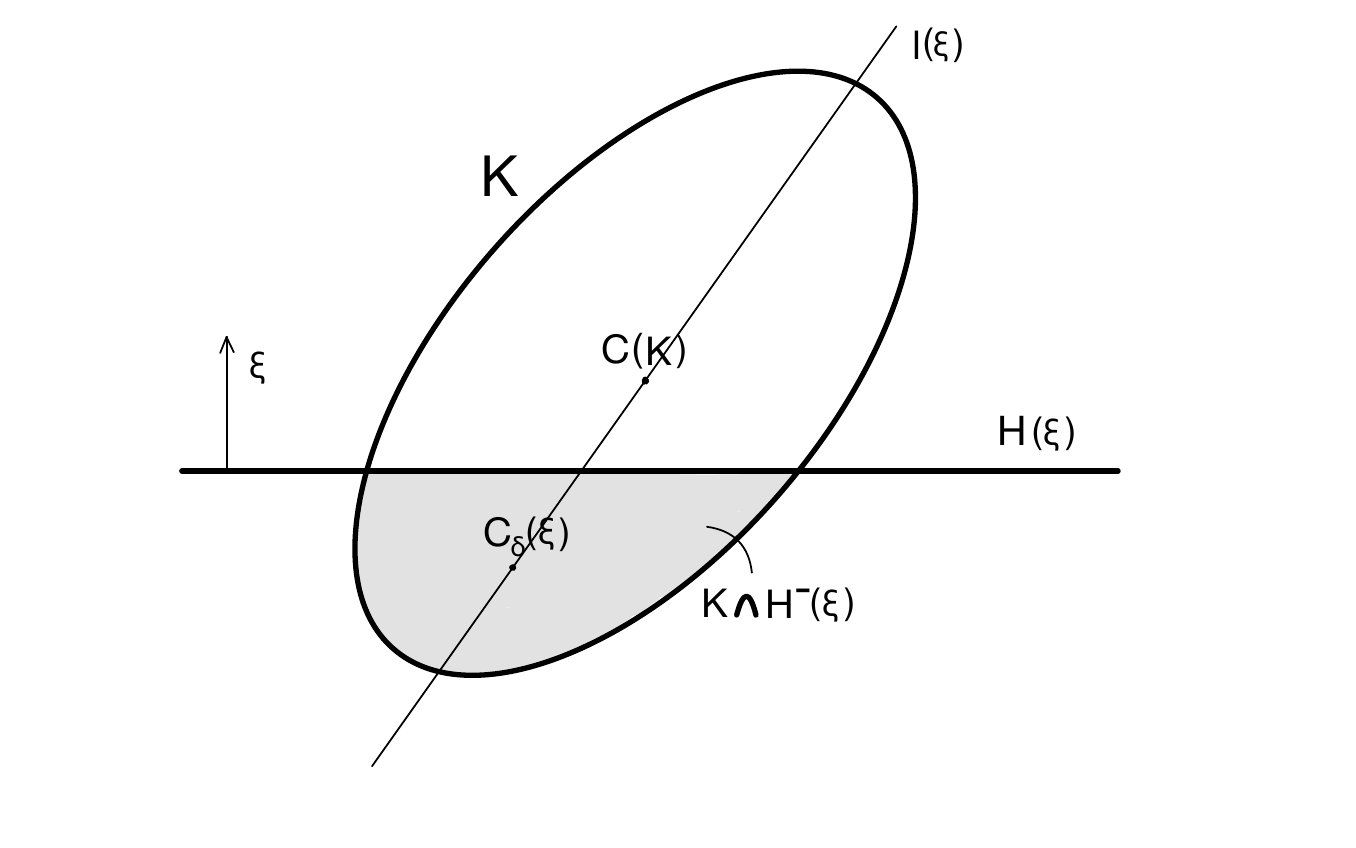} 
	\caption{A body $K$ and its submerged part  $K\cap H^-(\xi)$}
	\label{width1}
\end{figure}

We recall several well-known facts and definitions. 
Let $d\ge 3$, let $K\subset {\mathbb R^d}$ be a convex body 
and let $\delta\in (0, \textrm{vol}_d(K))$ be fixed. Given a direction $\xi\in S^{d-1}$ and $t=t(\xi)\in{\mathbb R}$, we call a hyperplane 
\begin{equation}\label{p2}
H(\xi)=H_t(\xi)=\{p\in{\mathbb R^{d}}:\,p\cdot\xi=t\},
\end{equation}
the {\it cutting hyperplane} of $K$ in the direction $\xi$,  if it cuts out of $K$ the given volume $\delta$, 
i.e., if 
\begin{equation}\label{fubu}
\textrm{vol}_d(K\cap H^-(\xi))=\delta,\qquad H^-(\xi)=\{p\in{\mathbb R^{d}}:\,p\cdot\xi\le t(\xi)\},
\end{equation}
(see Figure \ref{width1}).

Now we  recall the notions of {\it floating in equilibrium} and the {\it  surface of centers},  \cite{DVP}, \cite{Zh}.
\bd\label{efb}
Let   $\xi\in S^{d-1}$  and let ${\mathcal C}(\xi)={\mathcal C}_{\delta}(\xi)$ be the center of mass of the submerged part $K\cap H^-(\xi)$ satisfying (\ref{fubu}). 
A convex body $K$ {\it floats in equilibrium in the direction $\xi\in S^{d-1}$ at the level $\delta$} if (\ref{fubu}) holds and  the line $l(\xi)$ connecting  ${\mathcal C}(K)$ with ${\mathcal C}_{\delta}(\xi)$ is orthogonal to the ``free water surface" $H(\xi)$, i.e., the line $l(\xi)$ is ``vertical" $\textnormal{(}$parallel to $\xi$, see Figure \ref{width1}$\textnormal{)}$. We say that $K$ floats in equilibrium in every orientation at the level $\delta$ if $l(\xi)$ is parallel to $\xi$ for every $\xi\in S^{d-1}$.
\ed
\bd\label{scn}
Let $K$ be a convex body, let   $\xi\in S^{d-1}$  and let ${\mathcal C}(\xi)={\mathcal C}_{\delta}(\xi)$ be the center of mass of the submerged part $K\cap H^-(\xi)$ satisfying (\ref{fubu}). The geometric locus   $\{{\mathcal C}_{\delta}(\xi):\,\xi\in S^{d-1}\}$  is called the {\it surface of centers} ${\mathcal S}={\mathcal S}_{\delta}$ or the {\it  surface of buoyancy}  $\textnormal{(}$see Figure \ref{width2}$\textnormal{)}$. 
\ed

One can show, see Theorem \ref{D1} below, that the surface of centers is a  boundary of a strictly convex body.

\br\label{ew}
It was recently proved in \cite{HSW} that the surface of centers ${\mathcal S}$ is $C^{k+1}$-smooth, provided $K$ is of class $C^k$, $k\ge 0$.
In particular,  if $K$ is an arbitrary  convex body, then ${\mathcal S}$ is $C^1$-smooth.
\er

The following  result  is well-known,  see   
\cite[page 203]{G},   \cite[Section 2.1]{V1} and 
\cite[Corollary 2.4]{HSW}.
In the next section  we give a different proof.

\bl\label{tr}
Let  $d\ge 2$, let $K$ be a convex body  and let $\delta\in (0, \textnormal{vol}_d(K))$.
If $K$
floats in equilibrium in every orientation at the  level $\delta$, then   the surface of centers ${\mathcal S}$ is a sphere. Conversely, if ${\mathcal S}$ is a sphere centered at ${\mathcal C}(K)$, then $K$ floats in equilibrium in every orientation.
\el

It is known that the condition of ${\mathcal S}$ being  centered at ${\mathcal C}(K)$ is  satisfied for 
$\delta=\frac{\textnormal{vol}(K)}{2}$ (${\mathcal C}(K)$ is an arithmetic average of ${\mathcal C}(K\cap H^+(\xi))$ and ${\mathcal C}(K\cap H^-(\xi))$ for every $\xi\in S^{d-1}$),  and for any 
$\delta\in (0, \textrm{vol}_d(K))$, provided $K$ is  centrally-symmetric.

Now we recall the notion of a {\it floating body}.
A floating body $K_{[\delta]}$ of  $K$ was introduced by C. Dupin in 1822, \cite{D}.
\bd\label{dupa}
A non-empty convex set $K_{[\delta]}$ is the Dupin floating body of $K$ if each supporting plane of $K_{[\delta]}$ cuts off a set of volume $\delta\in (0, \textnormal{vol}_d(K))$ from $K$.
\ed
We remark that $K_{[\delta]}$ does not necessarily  exist for every convex $K$, see \cite{L} or \cite[Chapter 5]{NSW}, but
if $K$ has a sufficiently smooth boundary and $\delta>0$ is  small enough, then $K_{[\delta]}$ exists, \cite[Satz 2]{L}.

The notion of a {\it convex floating body}  was introduced independently in \cite{BL} and  \cite{SW1}. 
\bd\label{blsw}
A body $K_{\delta}$ is called the convex 
floating body of $K$, provided
$$
K_{\delta}=\bigcap\limits_{\{\xi\in S^{d-1}\}}H^+(\xi),\qquad H^+(\xi)=\{p\in{\mathbb R^d}:\,p\cdot\xi\ge t(\xi)\}.
$$
\ed
If $K_{[\delta]}$ exists, then $K_{[\delta]}=K_{\delta}$; 
$K_{\delta}$ is allowed to be an empty set, \cite{SW1}.
It was proved in \cite[Theorem 3, page 334]{MR} that $K_{[\delta]}=K_{\delta}$ for any 
$0<\delta\le\frac{\textrm{vol}_d(K_{\delta})}{2}$,
provided $K$ is centrally-symmetric. 
It was also shown in  \cite{MR}  that the boundary of $K_{\delta}$ is $C^2$-smooth, provided the boundary of $K$ is $C^1$-smooth and 
for every $x$ on the boundary of $K$  there is a unique supporting hyperplane
of $K$  through $x$.

Let $K$ float in equilibrium in every orientation for some   $\delta\in (0,\textrm{vol}_d(K))$, $\delta\neq \frac{\textrm{vol}_d(K)}{2}$.  It is not clear if the  additional  condition 
$K_{[{\delta}]}=K_{\delta}$
yields an affirmative answer to Ulam's Problem 19.

\subsection{Theorems of Dupin}
The solution of the problem  of finding the directions in which the given convex body floats in equilibrium is contained in the following  three results, proved by Dupin, (cf. \cite[pgs. 658-660]{Zh}  and \cite{Da} for $d=2$, and \cite[pgs. 287-288]{DVP} for $d=3$; see also \cite{G}). 
For convenience of the reader,  in this section we formulate these theorems for all $d\ge 3$ and include  sketches of the proofs in the next section.

\begin{figure}[h]
	\centering
	\includegraphics[height=1.5in]{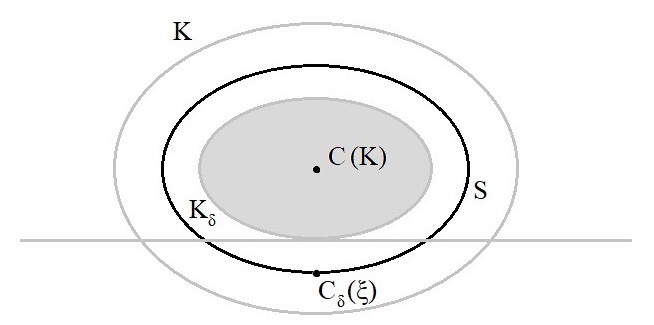} 
	\caption{Floating body $K_{\delta}$ and  surface of centers ${\mathcal S}$}
	\label{width2}
\end{figure}

Let $\xi\in S^{d-1}$ and let ${\mathcal H}(\xi)$ be a tangent hyperplane to ${\mathcal S}$ at ${\mathcal C}(\xi)$ which is the center of mass of $K\cap H^-(\xi)$, see Remark \ref{ew}.
The First Theorem of Dupin reads as follows.

\bt\label{D1}
Let $d\ge 2$, $K\subset {\mathbb R^d}$ be convex,  and let $\delta\in (0,\textnormal{vol}_d(K))$.  If   $H(\xi)$,  $\xi\in S^{d-1}$,  is a cutting hyperplane, then  ${\mathcal H}(\xi)$
 is parallel to 
 $H(\xi)$. Moreover, 
the bounded set $L({\mathcal S})$ with boundary ${\mathcal S}$ is a strictly convex body.
\et

The Second Theorem of Dupin is
\bt\label{D2}
Let $d\ge 2$, $K\subset {\mathbb R^d}$ be convex,  and let $\delta\in (0,\textnormal{vol}_d(K))$.
 Assume that $H(\xi)$,  $\xi\in S^{d-1}$, is a cutting hyperplane   and 
$\{H_n\}_{n=1}^{\infty}$, $H_n=H(\xi_n)$,  is any sequence of cutting hyperplanes converging to $H(\xi)$ as $\xi_n\to\xi$ for $n\to\infty$ 
and such that the limit $\lim\limits_{n\to\infty}H(\xi)\cap H(\xi_n)$ exists.
Then the  $(d-2)$-dimensional plane $\Pi=\lim\limits_{n\to\infty}H(\xi)\cap H(\xi_n)$  passes through the center of mass of $K\cap H(\xi)$. 
\et

\begin{figure}[h]
	\centering
	\includegraphics[height=3in]{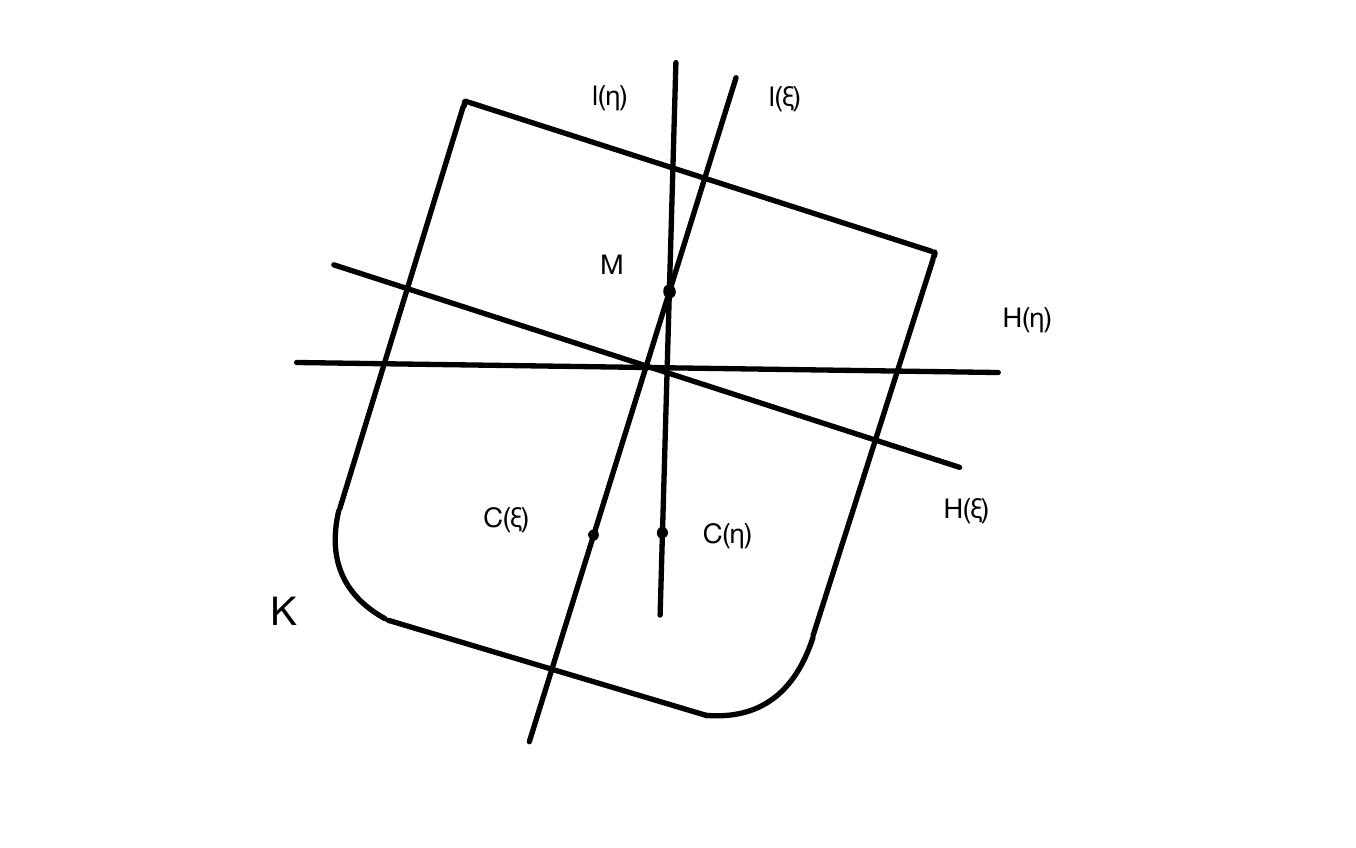} 
	\caption{The metacenter  $M=l(\xi)\cap l(\eta)$ of $K$}
	\label{meta}
\end{figure}

In order to formulate the third Theorem of Dupin in the case $d\ge 3$, we recall   the notions of a {\it metacenter} \cite[page 284]{DVP} and of a {\it moment of inertia} \cite[page 553]{Zh}.

To define the metacenter {\it heuristically}, assume 
 that  a body $K\subset {\mathbb R^3}$ is ``cylindrical".
In  naval architecture, \cite{Tu}, 
a ship floating originally at a horizontal waterline $H(\xi)\subset E$  is rotated through a small angle by an external force and then floats at waterline $H(\eta)\subset E$ (it is assumed that $H(\xi)$ and $H(\eta)$ intersect at the center of mass of $K$). Then the point $M=l(\xi)\cap l(\eta)$ is the metacenter, where
$l(\xi)$ is the line parallel to $\xi$ passing through 
the old center of boyancy ${\mathcal C}(\xi)$ and 
$l(\eta)$ is the line parallel to $\eta$ passing through 
the new center of boyancy ${\mathcal C}(\eta)$, see Figure \ref{meta}.

Now we recall a rigorous definition, \cite[pgs. 284, 285]{DVP}. 
\bd\label{MGF}
Let ${\mathcal S}$ be the surface of centers and 
let ${\mathcal C}$ be a point on ${\mathcal S}$ at which the normal curvatures exist.
Assume  that ${\mathcal C}$ belongs to 
some curve $\gamma\subset{\mathcal S}$ with the tangent $\zeta$ at ${\mathcal C}$. Take ${\mathcal C}'\in\gamma$ close to ${\mathcal C}$ and consider the normal lines $l_{\mathcal C}$, $l_{{\mathcal C}'}$,  to 
${\mathcal S}$ at ${\mathcal C}$ and ${\mathcal C}'$. If $\mu\mu'$ is a shortest distance between these lines, $\mu\in l_{\mathcal C}$, $\mu'\in l_{{\mathcal C}'}$, then the limiting position of the end $\mu$ of the segment 
$[\mu,\mu']$, when ${\mathcal C}'$ tends to ${\mathcal C}$,
is the metacenter $M_{{\mathcal C}}(\zeta)$ related to ${\mathcal C}$ in the tangential direction $\zeta$. 
\ed

Let ${\mathcal S}$ be $C^2$-smooth. One can assume without loss of generality that 
the tangent hyperplane ${\mathcal H}$ to ${\mathcal S}$ at ${\mathcal C}$ is horizontal, i.e., ${\mathcal H}$ is the $x_1\dots x_{d-1}$-hyperplane and that ${\mathcal C}$ is the origin. Then, choosing properly the directions of the axes in 
${\mathcal H}$ one can assume that the equation of ${\mathcal S}$ in a small neighborhood of ${\mathcal C}$ is
\begin{equation}\label{DVFrench1}
2x_d=k_1x_1^2+\dots +k_{d-1}x_{d-1}^2+o(x_1^2,\dots, x_{d-1}^2),
\end{equation}
where $k_j$, $j=1,\dots, d-1$, are  some non-negative coefficients, $k_1\le k_2\le$ $\dots\le k_{d-1}$.
\bl\label{DVPkr2}
The $x_d$-coordinate of $M_{{\mathcal C}}(\zeta)$  is
\begin{equation}\label{DVPkr1}
{\mathcal C}\mu=\frac{k_1\zeta_1^2+\dots+k_{d-1}\zeta_{d-1}^2}{k_1^2\zeta_1^2+\dots+k_{d-1}^2\zeta_{d-1}^2}, \qquad \textrm{where}\quad\zeta=(\zeta_1,\dots,\zeta_{d-1})\in S^{d-2}.
\end{equation}
\el
This formula is proved in \cite[page 285]{DVP} for $d=3$,  the general case can be shown similarly. For convenience of the reader we prove (\ref{DVPkr1}) in Appendix.
\br\label{Fkrut11}
We see that $\frac{1}{k_{d-1}}  \le {\mathcal C}\mu\le \frac{1}{k_1} $ and that ${\mathcal C}\mu$ is equal to one of $\frac{1}{k_{j}}$,
$j=1,\dots,d-1$, provided $\zeta$ is one of the corresponding principal directions of ${\mathcal S}$ at ${\mathcal C}$.
\er

We refer the reader to \cite[pgs. 103-106]{Sch2} and  \cite[pgs. 82-89]{T}  for the definition of the principal directions and  the normal  curvatures. 
Alexandrov proved that if $M$ is a  convex body and $G(\xi)$ is its supporting hyperplane,  then the normal curvatures exist  at $M\cap G(\xi)$ for almost every $\xi\in S^{d-1}$, \cite{BF}, \cite{Al}, \cite{H}. Hence, for an arbitrary convex body the metacenter is defined for almost every $\xi\in S^{d-1}$.

Now we define the moment of inertia.
Let $d\ge 3$, let $\delta\in (0,\frac{\textrm{vol}_d(K)}{2})$, and  let  $\xi\in S^{d-1}$ be any direction. Consider a convex body $K$ and the hyperplane  $H(\xi)$ defined by (\ref{p2}) such that   (\ref{fubu}) holds.
Choose  any $(d-2)$-dimensional plane $\Pi\subset H(\xi)$ passing through the center of mass ${\mathcal C}(K\cap H(\xi))$ and let 
$\eta_1,\dots,\eta_{d-2}, \eta_{d-1}$ be an orthonormal basis of $\xi^{\perp}=\{p\in{\mathbb R^d}:\,p\cdot\xi=0  \}$ such that 
\begin{equation}\label{base}
\Pi={\mathcal C}(K\cap H(\xi))+\textrm{span}(\eta_1,\dots,\eta_{d-2}),\quad H(\xi)={\mathcal C}(K\cap H(\xi))+\xi^{\perp}. 
\end{equation}

\bd\label{miab}
The moment of inertia $I_{K\cap H(\xi)}(\Pi)$ of $K\cap H(\xi)$ with respect to  $\Pi$  is calculated by summing $\textnormal{dist}(\Pi, v)^2$  for every ``particle" $v$ in the set $K\cap H(\xi)$, where $\textnormal{dist}(\Pi,v)=\min\limits_{\{x\in \Pi\} }|v-x|$, 
$\textnormal{(}$see Figure 4$\textnormal{)}$,
 i.e., 
\begin{equation}\label{moment}
I_{K\cap H(\xi)}(\Pi)=\int\limits_{K\cap H(\xi)}\textnormal{dist}(\Pi,v)^2dv=
\int\limits_{K\cap H(\xi)-{\mathcal C}(K\cap H(\xi))}(u\cdot\eta_{d-1})^2\,du.
\end{equation}
\ed

\begin{figure}[h]
	\centering
	\includegraphics[height=2.8in]{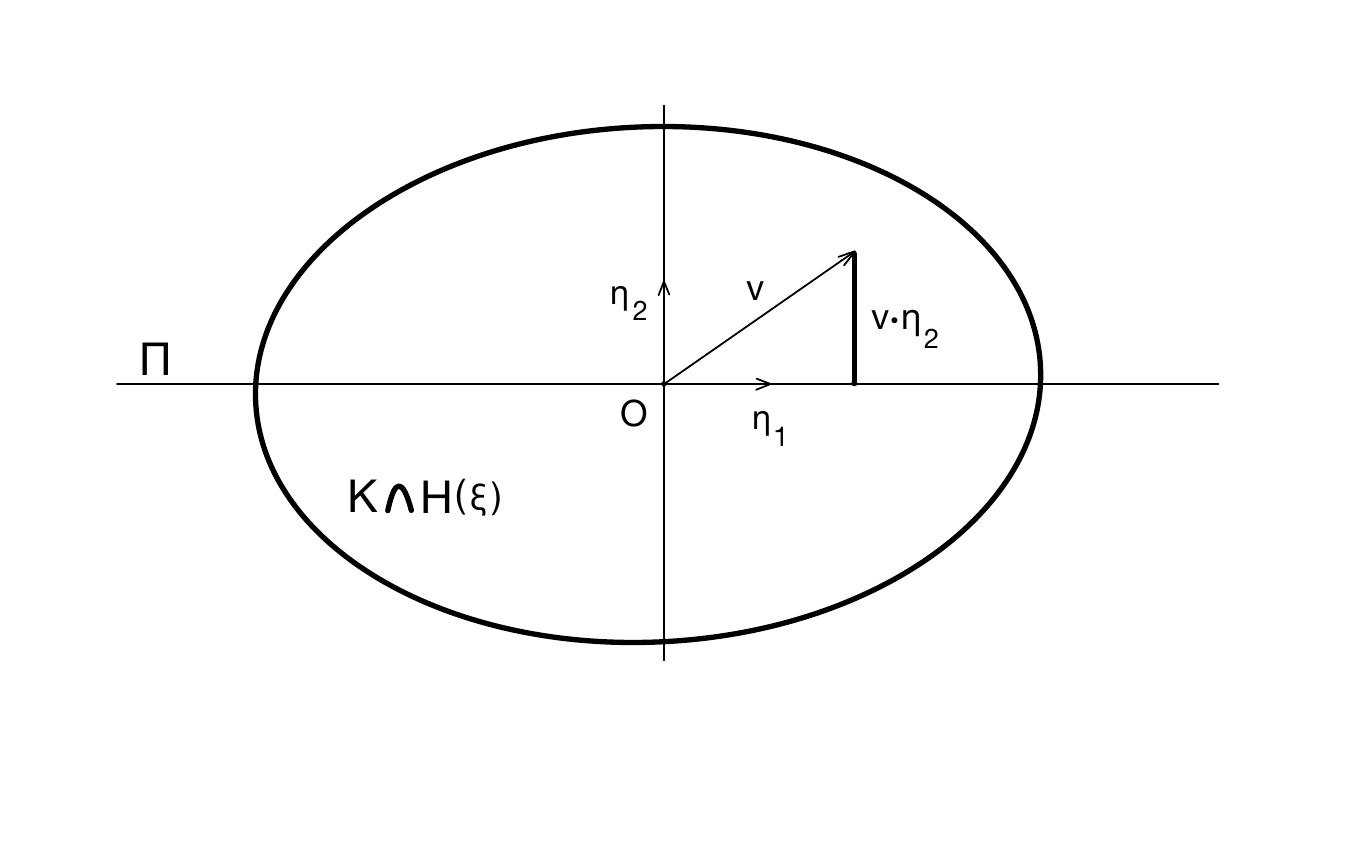} 
	\caption{Two-dimensional body $K\cap H(\xi)$ with  center of mass  at the origin, and a line $\Pi$ parallel to $\eta_1$; we have $\textrm{dist}(\Pi,v)^2=|v|^2-(v\cdot \eta_1)^2=(v\cdot \eta_2)^2$.}
	\label{width4}
\end{figure}

The Third  Theorem of Dupin reads as follows (cf. \cite{DVP}, page 288). 
\bt\label{D3}
Let $d\ge 3$, let $K\subset {\mathbb R^d}$ be a convex body  and let $\delta\in (0,\textnormal{vol}_d(K))$.
If  $H(\xi)$, $\xi\in S^{d-1}$, is a cutting  hyperplane  and  ${\mathcal C}={\mathcal C}(\xi) \in {\mathcal S}$ is the corresponding center of mass at which the normal curvatures of ${\mathcal S}$ exist in all directions and if a sequence of cutting hyperplanes $\{H_n\}_{n=1}^{\infty}$, $H_n=H(\xi_n)$,  converging to $H(\xi)$ as $n\to\infty$, is such that the limit  
$\lim\limits_{n\to\infty}H(\xi)\cap H(\xi_n)$
exists,
then for the corresponding sequence of the centers of mass $\{{\mathcal C}_n\}_{n=1}^{\infty}$, ${\mathcal C}_n={\mathcal C}(\xi_n)$,  ${\mathcal C}=\lim\limits_{n\to\infty}{\mathcal C}_n$,
one has
$$
{\mathcal R}_{{\mathcal C}(\xi)}(\zeta):=\textrm{dist}({\mathcal C}(\xi), M_{{\mathcal C}(\xi)}(\zeta))=\frac{1}{\delta} I_{K\cap H(\xi)}(\Pi),
$$
where  $\zeta=\lim\limits_{n\to\infty}\frac{{\mathcal C}{\mathcal C}_n}{|{\mathcal C}{\mathcal C}_n|}$ and $I_{K\cap H(\xi)}(\Pi)$ is the moment of inertia of 
$K\cap H(\xi)$ with respect to   the  $(d-2)$-dimensional plane $\Pi=\lim\limits_{n\to\infty}H(\xi)\cap H(\xi_n)$.
\et

If the reader does not want to deal with subtleties related to 
the almost everywhere existence of tangent hyperplanes or normal curvatures 
for general convex bodies, \cite{BF}, \cite{Al}, \cite{H}, one can assume from now on that   $K$ is  $C^1$. In this case,   ${\mathcal S}$ is $C^2$-smooth, \cite{HSW},  and Theorem \ref{D3} holds for every $\xi\in S^{d-1}$.

The following theorem can be found in \cite[page 23]{Da} and \cite{A} in the case when $K$ has $C^1$-smooth boundary. It  is  the Third Theorem of Dupin for $d=2$.
\bt\label{Dakr}
Let  $K\subset {\mathbb R^2}$ be convex  and  let  $\delta\in (0,\textnormal{area}(K))$. Then
$$
R(\xi)=\frac{\textnormal{length}^3(K\cap H(\xi))}{12\, \textnormal{area}(K\cap H^-(\xi))} \qquad\textrm{for almost every}\quad \xi\in S^1,
$$
where  $H(\xi)$  and $H^-(\xi)$ are defined by   (\ref{p2}) and (\ref{fubu}), and
$R(\xi)$ is the radius of curvature of ${\mathcal S}$ at the point of tangency ${\mathcal S}\cap {\mathcal H}(\xi)$.
\et

\section{Proofs of Theorems of Dupin}\label{DP}

\subsection{Proof of Theorem \ref{D1}} Rotating and translating if necessary  we can assume that  $\xi$ is such that $H(\xi)$ is ``horizontal", i.e., $H(\xi)=e_d^{\perp}$. Let $\eta\in S^{d-1}$, $\eta\neq \xi$ and let ${\mathcal H}(\xi)$ be a hyperplane parallel to $H(\xi)$ and passing through  ${\mathcal C}_{\delta}(\xi)$.  We claim that  ${\mathcal C}_{\delta}(\eta)$ is ``above" 
${\mathcal H}(\xi)$, i.e., 
$x_d({\mathcal C}_{\delta}(\xi))<x_d({\mathcal C}_{\delta}(\eta))$.
Since  $x_d>0$
$\forall x\in (K\cap H^-(\eta))\setminus  (K\cap H^-(\xi))$ but $x_d\le 0$ 
$\forall x\in (K\cap H^-(\xi))\setminus  (K\cap H^-(\eta)$, we have
$$
x_d({\mathcal C}_{\delta}(\xi))=\frac{1}{\delta}\Big( \int\limits_{(K\cap H^-(\xi))\setminus  (K\cap H^-(\eta))}x_d dx+\int\limits_{K\cap H^-(\eta)\cap H^-(\xi)}x_d dx\Big)<
$$
$$
\frac{1}{\delta}\Big( \int\limits_{(K\cap H^-(\eta))\setminus  (K\cap H^-(\xi))}x_d dx+\int\limits_{K\cap H^-(\eta)\cap H^-(\xi)}x_d dx\Big)=\,x_d({\mathcal C}({\mathcal C}_{\delta}(\eta))
$$
and  the claim is proved.
Thus, for any $\xi\in S^{d-1}$  we have
${\mathcal S}\subset {\mathcal H}^+(\xi)$,
${\mathcal S}\cap {\mathcal H}(\xi)={\mathcal C}_{\delta}(\xi)$ and $\min\limits_{\{\xi\in S^{d-1}\}}|{\mathcal C}(K)-{\mathcal C}_{\delta}(\xi)|>0$. We conclude that
$L({\mathcal S})=\bigcap\limits_{\{\xi\in S^{d-1}\}}{\mathcal H}^+(\xi)$ is a strictly convex body.
$\qquad\qquad\qquad\qquad\qquad\qquad\qquad \,\,\,\square $

\subsection{Proof of Theorem \ref{D2}}

Rotating and translating if necessary, assume that $H(\xi)$ is ``horizontal", i.e., $H(\xi)=e_d^{\perp}$. Take $n$ large enough and consider the $(d-2)$-dimensional plane $\Pi_n=H(\xi)\cap H(\xi_n)$. Introduce the ``moving" coordinates $(x_1,x_2,\dots, x_{d-1}, x_d)$ so that  $\Pi_n$ is the 
$(x_2,\dots, x_{d-1})$-plane.

\begin{figure}[ht]
	\includegraphics[width=360pt]{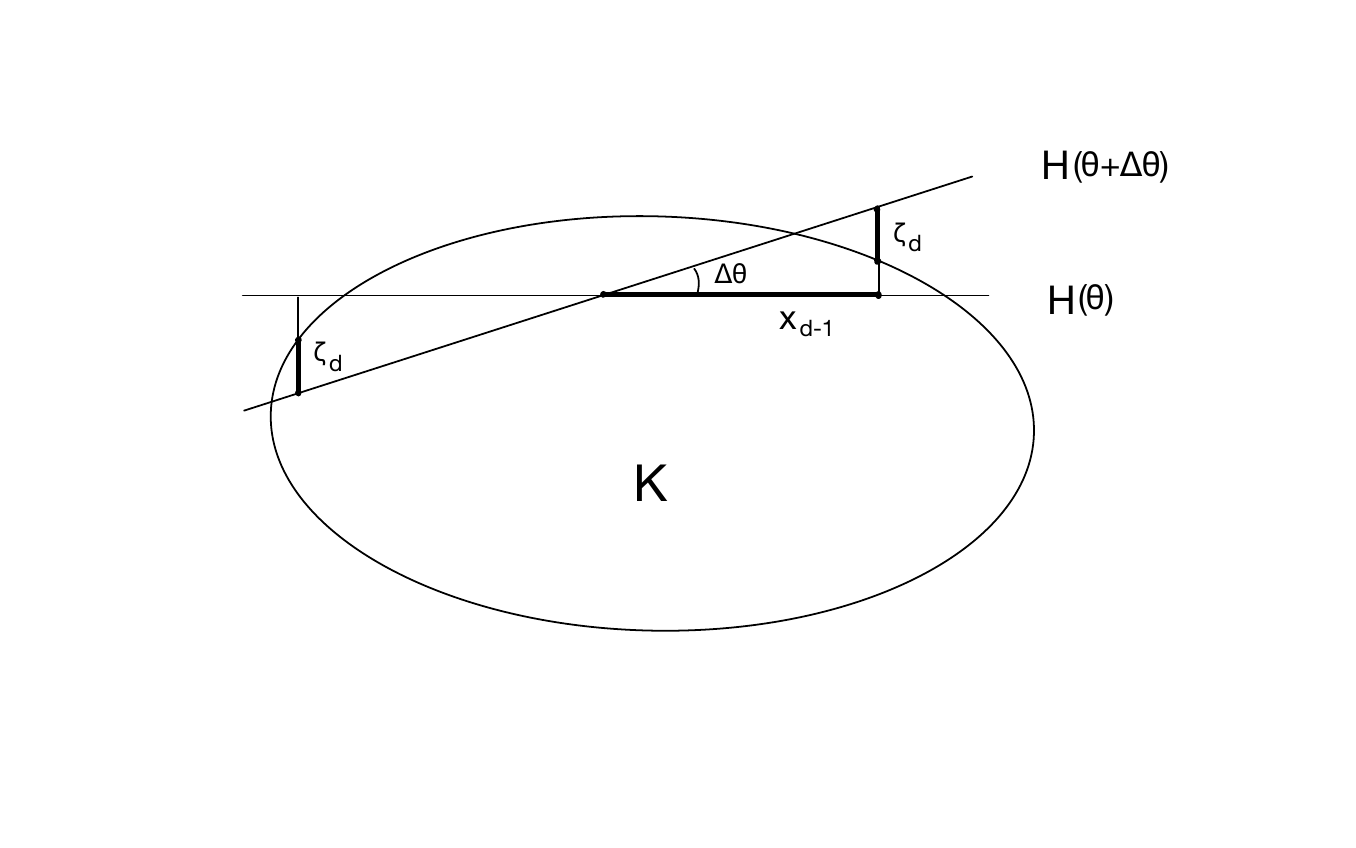}
	\caption{The  function $\zeta_d$.}
	\label{Fig7}
\end{figure}

Denote by $A\triangle B$ the symmetric difference of two sets $A$ and $B$, i.e., 
$A\triangle B=(A\setminus B)\cup (B\setminus A)$, and 
let
$\Lambda_n= (K\cap H(\xi))\triangle P_{H(\xi)} (K\cap H(\xi_n))
$, where
$P_{H(\xi)}$ is the orthogonal projection onto $H(\xi)$. Then,
\begin{equation}\label{kr1}
\varDelta V=\textnormal{vol}_d(K\cap H^-(\xi))-\textnormal{vol}_d(K\cap H^-(\xi_n))=
\end{equation}
$$
\int\limits_{K\cap H(\xi)}x_1\tan \varepsilon_n \,dx-
\int\limits_{\Lambda_n}\zeta_d \,dx=0,
$$
where  $x_1=x_1(\xi,\xi_n)$ and 
$\zeta_d=\zeta_d(\xi,\xi_n)$ is an error  of $x_d=x_1\tan\varepsilon_n$ in $\Lambda_n$ which is obtained  during the computation of $\varDelta V$ using the first integral above (see Figure \ref{Fig7}; observe that 
$H(\xi)\cap H(\xi_n)\cap \textrm{int}K\neq\emptyset$
(see  \cite[p. 116]{O} or \cite[Appendix A]{R2})).               
To see (\ref{kr1}),   consider on $e_d^{\perp}$ an infinitesimally small element of the $(d-1)$-dimensional volume $dx$ as a base of an infinitesimally small prism ``between"
$H(\xi)$ and  $H(\xi_n)$ of ``height" $\tan\varepsilon_n |x_1|$, where  $\varepsilon_n$ is a small angle between 
$H(\xi)$ and  $H(\xi_n)$. The $d$-dimensional volume of the prism  is $\tan\varepsilon_n|x_1|dx$.  Summing up the  volumes of the corresponding prisms we obtain (\ref{kr1}).

By (\ref{kr1}), we have
$$
x_1({\mathcal C}(K\cap H(\xi))=\frac{\int\limits_{K\cap H(\xi)}x_1 \,\,\,dx}{\textrm{vol}_{d-1}(K\cap H(\xi))}=\frac{\int\limits_{\Lambda_n}\zeta_d \,\,\,dx}{\textrm{vol}_{d-1}(K\cap H(\xi))\tan \varepsilon_n}.
$$
Since $\textrm{vol}_{d-1}(\Lambda_n)\to 0$ as $n\to\infty$ (see \cite[p. 116]{O} or \cite[Appendix A]{R2}), and since $|\zeta_d|\le D\tan\varepsilon_n$, where $D$ is the diameter of $K$, 
we obtain
$$
|x_1({\mathcal C}(K\cap H(\theta)))|\le \frac{D\tan \varepsilon_n\,\,\textrm{vol}_{d-1}(\Lambda_n)}{\textrm{vol}_{d-1}(K\cap H(\xi))\tan \varepsilon_n}\to 0
$$
as $n\to \infty$.
We see that the $(d-2)$-dimensional plane 
$H(\xi)\cap H(\xi_n)$ tends, as $n\to \infty$,  to a limiting position  $\Pi$  that passes through the center of mass of $K\cap H(\xi)$. 
$\qquad\qquad\qquad\qquad \qquad\qquad\qquad\qquad\qquad \qquad\qquad\qquad\qquad\, \square $

\begin{figure}[h]
	\centering
	\includegraphics[height=3.5in]{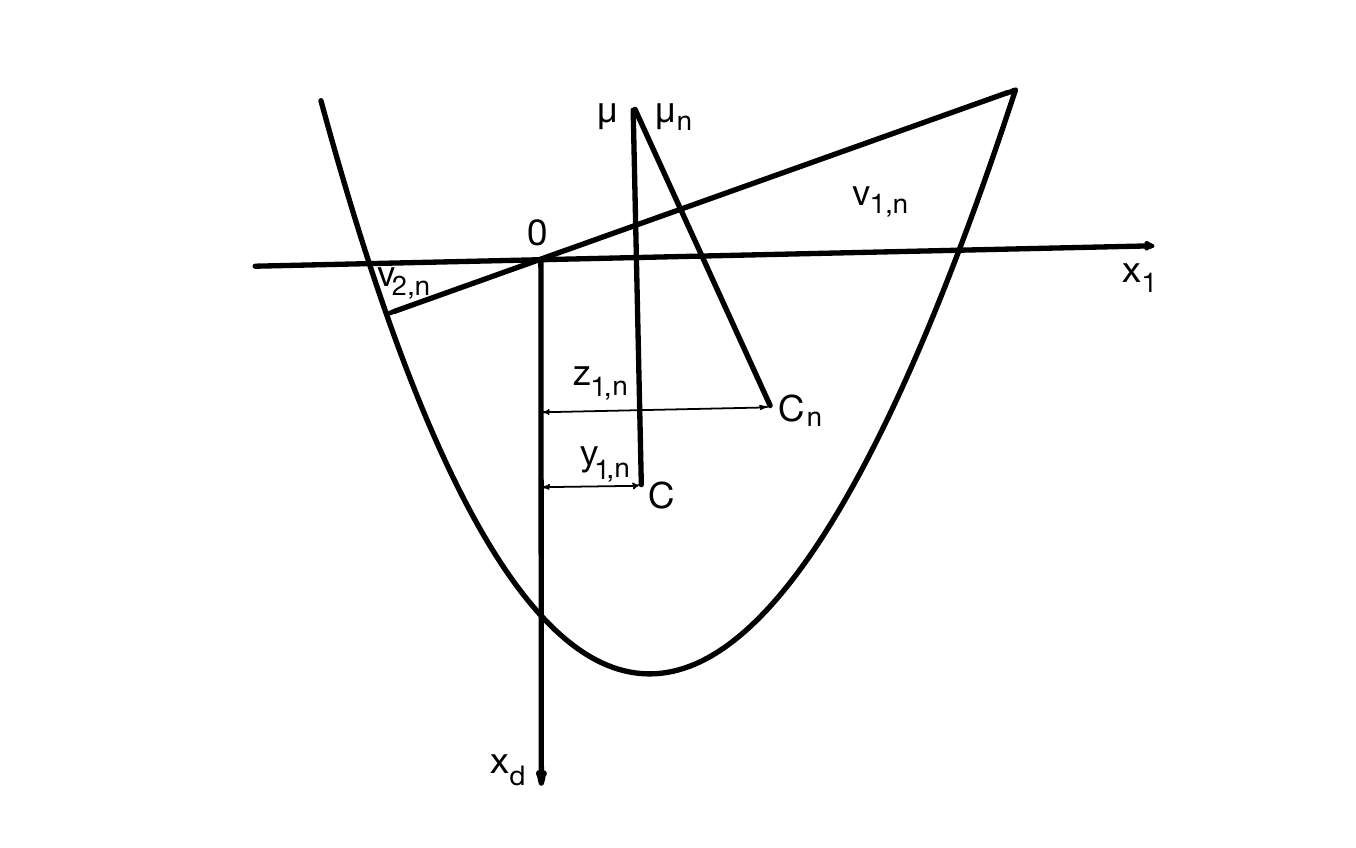} 
	\caption{The normals ${\mathcal C}\mu$ and ${\mathcal C}_n\mu_n$ to the surface of centers}
	\label{VPbeauty}
\end{figure}

\subsection{Proof of Theorem \ref{D3}}

As in the previous proofs, we assume that $H(\xi)=e_d^{\perp}$. We take $n$ large enough and put $\Pi_n=H(\xi)\cap H(\xi_n)$.  As above we introduce the ``moving" coordinates $(x_1,x_2,\dots, x_{d-1}, x_d)$ so that the $(d-2)$-dimensional plane $\Pi_n$ is the 
$(x_2,\dots, x_{d-1})$-plane. Denote by $v_{1,n}$ and $v_{2,n}$ the $d$-dimensional bodies with the $x_1$-coordinates having opposite 
signs,
$$
v_{1,n}=(K\cap H^-(\xi_n))\setminus  (K\cap H^-(\xi)),\quad v_{2,n}=(K\cap H^-(\xi))\setminus  (K\cap H^-(\xi_n)),
$$
 and let $y_{1,n}$, $z_{1,n}$ be the $x_1$-coordinates of  ${\mathcal C}={\mathcal C}_{\delta}(\xi)$ and ${\mathcal C}_n={\mathcal C}_{\delta}(\xi_n)$, see Figure \ref{VPbeauty} (cf. Figure 59, page 289 from \cite{DVP}). Then
$$
\delta y_{1,n}=\int\limits_{ K\cap H^-(\xi)}x_1dx,\quad \delta z_{1,n}=\int\limits_{ K\cap H^-(\xi_n)}x_1dx,
$$
and looking at the difference, we have
$$
\delta (y_{1,n}-z_{1,n})=\int\limits_{v_{1,n}\cup\, v_{2,n}}|x_1|dx.
$$
Repeating the argument from the proof of Theorem \ref{D2}  showing that the  volumes $\textrm{vol}_d(v_{1,n})=\textrm{vol}_d(v_{2,n})$
  are  (up to $o(\varepsilon_n)$) the sums of  volumes $\varepsilon_n x_1dx$ of infinitesimal prisms, we obtain
  \begin{equation}\label{dvp33}
  \delta (z_{1,n}-y_{1,n}) =\tan\varepsilon_n \int\limits_{K\cap H(\xi)}x_1^2d\sigma_{d-1}(x)+o(\varepsilon_n)=
  \end{equation}
$$
\tan\varepsilon_n I_{K\cap H(\xi)}(\Pi_n)+o(\varepsilon_n).
$$
On the other hand, consider the normals ${\mathcal C}\mu$ and ${\mathcal C}_n\mu_n$ to ${\mathcal S}$ at the points ${\mathcal C}={\mathcal C}_{\delta}(\xi)$ and ${\mathcal C}_n={\mathcal C}_{\delta}(\xi_n)$.  The angle $\varepsilon_n$ between these normals
is equal to the one between the   hyperplanes 
$H(\xi)$ and $H(\xi_n)$. At the same time this is the angle between the $x_d$-axis and ${\mathcal C}_n\mu_n$. By  definition of the metacenter, the vector $\mu\mu_n$ is ``parallel" to $\Pi_n$, so $\mu$ and $\mu_n$ have the same $x_1$-coordinate;
it is the $x_1$-coordinate of the intersection of orthogonal projections of lines $\ell$, $\ell_n$, containing 
${\mathcal C}\mu$, ${\mathcal C}_n\mu_n$, onto the $x_1x_d$-plane.
We conclude  that  $z_{1,n}-y_{1,n}$ is the projection of ${\mathcal C}_n\mu_n$ onto the $x_1$-axis, i.e.,
$z_{1,n}-y_{1,n}=\sin\varepsilon_n |{\mathcal C}_n\mu_n|$.
Substituting this expression into (\ref{dvp33}) and passing to the limit as $n\to\infty$ we see that
$$
|{\mathcal C}\mu|=\lim\limits_{n\to\infty}|{\mathcal C}_n\mu_n|=\frac{I_{K\cap H(\xi)}(\Pi)}{\delta},
$$
which is the desired conclusion. $\qquad\qquad\qquad\qquad \qquad\qquad\qquad\qquad\qquad \square $

\section{Proofs of Lemma \ref{tr}, Theorems \ref{Fedja1},  \ref{CS}, and Corollaries \ref{Al}, \ref{equi}}\label{Au}

We start with the proof of Lemma \ref{tr} (cf. \cite{Gr}, \cite{Mo},  \cite[page 203]{G} and \cite[Corollary 2.4 and Proposition 2.2]{HSW}).

\bp
 At first we  prove the  converse statement. Using  the fact that all normals of the sphere intersect at its center and Theorem \ref{D1},  we see that for  every $\xi\in S^{d-1}$, the lines $\ell(\xi)$ passing through ${\mathcal C}(K)$ and ${\mathcal C}_{\delta}(\xi)$ are orthogonal to $H(\xi)$.

Now we prove the {\it if} part.
Let $\xi\in S^{d-1}$ and let $\ell(\xi)$ be a line  passing through ${\mathcal C}(K)$ and the center of mass ${\mathcal C}(\xi)$ of $K\cap H^-(\xi)$.  
By
Theorem \ref{D1}, ${\mathcal H}(\xi)$ is parallel to  $H(\xi)$.
Since $K$ floats in equilibrium in the direction $\xi$,  the line $\ell(\xi)$  is orthogonal to  $H(\xi)$. Since ${\mathcal H}(\xi)$ is parallel to  $H(\xi)$,  $\ell(\xi)$  is the normal line to
 ${\mathcal S}$ at ${\mathcal C}(\xi)$, and since  the body floats in equilibrium in all directions $\xi\in S^{d-1}$, 
  we know that 
 the lines $\ell(\xi)$  passing through  ${\mathcal C}(K)$
 are the normal lines  to  ${\mathcal S}$ for every  $\xi$;  
 we recall that   ${\mathcal S}$  is $C^1$-smooth, \cite{HSW}. Consider any two-dimensional plane $\Pi$ passing through ${\mathcal C}(K)$.  Parametrizing   the plane curve
 $ {\mathcal S}\cap\Pi$ by the radius vector ${\mathbf r}$ going from ${\mathcal C} (K)$ to the corresponding 
 $ {\mathcal S}\cap l(\xi)$, we see that
 ${\mathbf r}$ is orthogonal to ${\mathbf r}'$, i.e., 
 ${\mathbf r}\cdot {\mathbf r}'=0$, $|{\mathbf r}|$ is constant, 
 and $ {\mathcal S}\cap \Pi$ is a circle.
 Since $\Pi$ was chosen arbitrarily, 
 applying \cite[Corollary 7.1.4, page 272]{Ga} to $L({\mathcal S})$ from Theorem \ref{D1}, we obtain that ${\mathcal S}$ is a sphere.
 This gives  the desired conclusion.
 \ep

\subsection{Proof of Corollary  \ref{Al}}
Let $\delta_n\to0$ and let ${\mathcal S}_n$ be the corresponding surfaces of centers, which are all spheres of the radii $r_n$, $r_n\to r$ as $n\to\infty$.
Since $d(K_{{\delta}_n}, K)\to 0$ as $n\to \infty$, and since $K_{{\delta}_n}\subset B_{r_n}^2(0)\subset K$, we have
$d(B_{r_n}^2(0), K)\to 0$ as $n\to \infty$. Hence, $K$ is the Euclidean ball $B^2_r(0)$. 
$\qquad\,\, \square$

\subsection{Proof of Theorem \ref{Dpr}}

It is a consequence of  Lemma \ref{tr} and Theorems of Dupin.
It will be convenient to reformulate Theorem \ref{Dpr} in terms of the radial function.

Given a direction $\xi\in S^{d-1}$ and a hyperplane  (\ref{p2}) for which  (\ref{fubu}) holds,
we will use  the notation $\rho_{K\cap H(\xi)}(w)$ for the radial function of the $(d-1)$-dimensional convex body $K\cap H(\xi)$  with respect to the center of mass ${\mathcal C}(K\cap H(\xi))$ in the direction $w\in S^{d-1}\cap \xi^{\perp}$, i.e., for
$$
\rho_{K\cap H(\xi),\,{\mathcal C}(K\cap H(\xi)) }(w)=\max\{\lambda>0:\,{\mathcal C}(K\cap H(\xi))+\lambda w\in (K\cap H(\xi))\}.
$$

\bt\label{Fedja1}
Let  $d\ge 3$, let $K$ be a convex body  and let $\delta\in (0, \textnormal{vol}_d(K))$.
If $K$ floats in equilibrium at the level $\delta$ in every orientation, then $\forall \xi\in S^{d-1}$ 
the cutting sections $K\cap H(\xi)$ have equal principal moments, i.e., 
we have
\begin{equation}\label{eq112}
\int\limits_{S^{d-1}\cap \xi^{\perp}}w_k^2\,\rho_{K\cap H(\xi)}^{d+1}(w)dw=(d+1)\delta  {\mathcal R},\quad k=1,2,\dots, d-1,
\end{equation}
\begin{equation}\label{eq112B}
\int\limits_{S^{d-1}\cap \xi^{\perp}}w_jw_k\,\rho_{K\cap H(\xi)}^{d+1}(w)dw=0,\quad 1\le k,j\le d-1, \quad j\neq k,
\end{equation}
where  ${\mathcal R}$ 
is the radius of the spherical surface of centers ${\mathcal S}$.

Conversely,  if ${\mathcal C}({\mathcal S})={\mathcal C}(K)$
and for every cutting hyperplane $H(\xi)$, $\xi\in S^{d-1}$, the cutting section $K\cap H(\xi)$ satisfies (\ref{fubu}), (\ref{eq112}) and (\ref{eq112B}) with some constant ${\mathcal R}$, then the body $K$ with  $C^1$-smooth boundary floats in equilibrium in every orientation 
 at the level $\delta$.
\et
\bp
 Let $d\ge 3$. 
Fix any $\xi\in S^{d-1}$ and a cutting hyperplane $H(\xi)$.
Let  $\Pi\subset H(\xi)$ be a $(d-2)$-dimensional plane  passing through ${\mathcal C}(K\cap H(\xi))$, 
let  $\Pi_n\subset H(\xi)$ be a sequence of $(d-2)$-dimensional planes converging  and parallel to $\Pi$ as $n\to\infty$, and let $H_n=H(\xi_n)$, $H_n\cap H(\xi)=\Pi_n$, be the corresponding cutting hyperplanes.
If  ${\mathcal C}_n={\mathcal C}(\xi_n)$ are the  centers of mass of $K\cap H_n^-$ converging to ${\mathcal C}={\mathcal C}(\xi)$ as $n\to\infty$, then, 
by
Theorem \ref{D3},
for $\zeta=\lim\limits_{n\to\infty}\frac{{\mathcal C}{\mathcal C}_n}{|{\mathcal C}{\mathcal C}_n|}$ we have
\begin{equation}\label{huhuh}
{\mathcal R}_{\mathcal C(\xi)}(\zeta)\stackrel{\text{for a.e \,$\xi$}}{=}\frac{1}{\delta}I_{K\cap H(\xi)}(\Pi).
\end{equation}

By Lemma \ref{tr} the surface of centers ${\mathcal S}$ is a sphere of  certain radius ${\mathcal R}$ centered at ${\mathcal C}(K)$.  Since the radii of the normal curvatures of  the sphere of radius ${\mathcal R}$ are equal to ${\mathcal R}$ at all points ${\mathcal C}\in{\mathcal S}$  in all directions and since $\Pi$ was chosen arbitrarily, 
by Remark \ref{Fkrut11}, we see that   the function in the right-hand side of (\ref{huhuh}) is constant for almost every $\xi\in S^{d-1}$ and for all $\Pi$.
Since the function 
$(\xi,\Pi)\to I_{K\cap H(\xi)}(\Pi)$
is continuous, 
the right-hand side of (\ref{huhuh}) is constant for  every $\xi\in S^{d-1}$ and for all $\Pi$.

Hence, using  (\ref{moment})  we  obtain that for all $\xi\in S^{d-2}$ one has
\begin{equation}\label{cognac}
\frac{1}{\delta}\int\limits_{K\cap H(\xi)-{\mathcal C}(K\cap H(\xi))}(v\cdot\eta_{d-1})^2\,dv={\mathcal R}\qquad\forall \eta_{d-1}\in S^{d-1}\cap \xi^{\perp},
\end{equation}
where  we recall that $\eta_1,\dots,\eta_{d-2}, \eta_{d-1}$ is the orthonormal basis of $\xi^{\perp}$ such that (\ref{base}) holds.
Passing to  polar coordinates in $H(\xi)$ with respect to ${\mathcal C}(K\cap H(\xi))$, we have
\begin{equation}\label{cognac22}
\int\limits_{K\cap H(\xi)-{\mathcal C}(K\cap H(\xi))}\!\!(v\cdot\eta_{d-1})^2dv=\int\limits_{S^{d-1}\cap \xi^{\perp}}\!\!dw\!\!\int\limits_0^{\rho_{K\cap H(\xi)}(w)}\!\!( r w\cdot \eta_{d-1})^2r^{d-2}dr=
\end{equation}
$$
\frac{1}{d+1}\int\limits_{S^{d-1}\cap \xi^{\perp}}(w\cdot\eta_{d-1})^2\rho^{d+1}_{K\cap H(\xi)}(w)dw,
\qquad\forall \eta_{d-1}\in S^{d-1}\cap \xi^{\perp}.
$$
This identity and  (\ref{cognac}) yield
\begin{equation}\label{vodochka1}
\int\limits_{S^{d-1}\cap \xi^{\perp}}(w\cdot\eta_{d-1})^2\rho^{d+1}_{K\cap H(\xi)}(w)dw=(d+1)\,\delta \,{\mathcal R},
\end{equation}
where the right-hand side is independent of $\eta_{d-1}\in S^{d-1}\cap \xi^{\perp}$. 
By choosing $\eta_{d-1}$ to be the standard coordinate vectors in $\xi^{\perp}$,
we obtain  (\ref{eq112}).  By taking $\eta_{d-1}=(0,\dots,\underbrace{\frac{\sqrt{2}}{2}}_j,0,\dots,0,\underbrace{\frac{\sqrt{2}}{2}}_k, 0,\dots, 0 )$ for different $1\le j,k\le d-1$, $j\neq k$, and using (\ref{eq112}) we obtain (\ref{eq112B}).
Since $\xi$ was arbitrary, the proof of the {\it if} part  is complete.

Now we prove the converse statement. Our goal is to show that the surface of centers is a sphere. 

We will show at first that 
for almost every $\xi\in S^{d-1}$ the points 
${\mathcal C}(\xi)={\mathcal S}\cap{\mathcal H}(\xi)$ are umbilical.
Let $\xi\in S^{d-1}$ be such that the normal curvatures at the corresponding point ${\mathcal C}(\xi)\in{\mathcal S}$ exist.
Assume  that (\ref{eq112}) and (\ref{eq112B}) are true.  We can also assume  that   $\Pi$ satisfies (\ref{base}).  
Then, expanding the expression $(w\cdot \eta_{d-1})^2$ by writing $w$ in the basis $\eta_1$, $\dots$, $\eta_{d-1}$
and 
using
the identities (\ref{cognac}) and (\ref{cognac22}), we see that (\ref{vodochka1}) holds with some constant ${\mathcal R}$ in 
the right-hand side, i.e., it is   independent of $\eta_{d-1}\in S^{d-1}\cap \xi^{\perp}$. 
Hence, using (\ref{moment}), (\ref{cognac}) and (\ref{cognac22}), we see that the right-hand side of 
(\ref{huhuh}) is independent of $\Pi$ and $\xi$.

Now let $\zeta$ be any unit principal direction in the hyperplane ${\mathcal H}(\xi)$ tangent to ${\mathcal S}$ at ${\mathcal C}(\xi)$, and let $\Pi$ be  a two-dimensional subspace  spanned by $\zeta$ and the normal  to ${\mathcal S}$ at ${\mathcal C}(\xi)$. Consider 
 a sequence of unit directions $\zeta_n$
 tangent to the two-dimensional curve ${\mathcal S}\cap\Pi$ at the corresponding points 
 ${\mathcal C}(\xi_n)\in({\mathcal S}\cap\Pi)$
 and such that $\zeta_n\to\zeta$, ${\mathcal C}(\xi_n)\to {\mathcal C}(\xi)$, as $n\to\infty$.
 If $\{H(\xi_n)\}_{n=1}^{\infty}$
 is
 a sequence of cutting hyperplanes $H(\xi_n)$ converging to $H(\xi)$ as $n\to\infty$ with  ${\mathcal C}(\xi_n)$ being the centers of mass of $K\cap H^-(\xi_n)$, applying
 Theorem \ref{D3} and passing to a subsequence if necessary to ensure the existence of $\lim\limits_{n\to\infty}H(\xi)\cap H(\xi_n)$,
 we see that the radii of the principal normal curvatures of 
 ${\mathcal S}$ at ${\mathcal C}(\xi)$ in the principal directions are  the same  and the value of the radii  is  independent of $\xi$ and $\zeta$
for almost every $\xi\in S^{d-1}$  
and for every principal direction $\zeta$ parallel to  ${\mathcal H}(\xi)$.

Thus, for almost every $\xi\in S^{d-1}$ the points 
${\mathcal C}(\xi)$ are umbilical.
We claim that ${\mathcal S}$ is a sphere.
Indeed,  recall that by Remark \ref{ew} the surface of centers is $C^2$.  Hence, by continuity, all the points on ${\mathcal S}$ are umbilical.  Using   \cite[Proposition 4, page 147]{DC} and  \cite[Corollary 7.1.4, page 272]{Ga}  we conclude that 
${\mathcal S}$ must be a $(d-1)$-dimensional sphere.
An application of Lemma \ref{tr} finishes the proof. \ep

\br\label{nun}
In the planar case an analogous result is a  consequence of Lemma \ref{tr} and Theorem 
\ref{Dakr}.
\er

\subsection{Proof of Corollary \ref{equi}}

The condition of the corollary reads as
\begin{equation}\label{radhop1}
\forall \xi\in S^{d-1},\qquad   \rho^{d+1}_{K\cap H(\xi)}(w)+ \rho^{d+1}_{K\cap H(\xi)}(-w)=c\qquad\forall w\in S^{d-1}\cap\xi^{\perp}.
\end{equation}
The result  follows  from the second part of Theorem \ref{Al} by writing $\rho_K^{d+1}$ as the sum of even and odd parts and
substituting the even part  from (\ref{radhop1}) into (\ref{eq112}) and  (\ref{eq112B}).$\qquad\qquad\qquad\qquad\qquad\qquad\qquad\qquad\qquad\qquad\qquad\qquad\qquad\quad\quad\,\,\,\,
\square$

\subsection{Proof of Theorem  \ref{CS}}

We recall that a measurable function $f:\, S^{d-1}\to {\mathbb R}$ is isotropic if the signed measure $fdx$  is isotropic,
i.e.,  its center of mass is at the
origin and the map
$$
S^{d-1}\ni y\quad \to\quad \int\limits_{S^{d-1}}(y\cdot w)^2f(w)dw
$$
is constant, \cite{MP}. The following result was obtained in  \cite{MRS}.

\bt\label{CrSa}
Let $f:\, S^{d-1}\to {\mathbb R}$ be a measurable, bounded a. e. and even function, $d\ge 3$. If
for almost every $\xi\in S^{d-1}$ the restriction $f|_{S^{d-1}\cap \xi^{\perp}}$ to  
$S^{d-1}\cap \xi^{\perp}$
is isotropic $\textnormal{(}$i.e. the
restriction of $f$ to almost every equator is isotropic$\textnormal{)}$, then $f$ is almost everywhere equal to a
constant.
\et

By the origin-symmetry,
the centers of mass of all cutting sections are equal to the center of  mass of $K$. Hence, we may apply
 Theorem \ref{Fedja1} to see that there exists a constant $c$ such that all second moments of the central sections $K\cap\xi^{\perp}$  are equal to $c$  for all $ \xi\in S^{d-1}$. The result follows from Theorem \ref{CrSa} with $f=\rho_K^{d+1}$. $\qquad\qquad\qquad\qquad\qquad\qquad\qquad\qquad\qquad\qquad\,\square$

\section{Appendix: proof of Lemma \ref{DVPkr2} from \cite[page 285]{DVP}}

Let $M$ be a point on $C^2$-smooth ${\mathcal S}$  and let $\gamma\subset {\mathcal S}$ be a curve passing through $M$. Let $M'\in\gamma$ be a point infinitesimally close to $M$. Consider two normal lines $M\Gamma$ and $M'N'$ to ${\mathcal S}$ at $M$ and $M'$ and let $\mu\mu'$ be the shortest distance between these normal lines. 
We can assume that the tangent hyperplane to ${\mathcal S}$ at $M$ is $e_d^{\perp}$ and that its boundary is locally described by (\ref{DVFrench1}).

Now  drop the terms of the orders higher than $2$.
We have $\frac{\partial x_d}{\partial x_j}=k_jx_j$ for $j=1,\dots,d-1$.
The normal line at $M'=M'(x_1,\dots,x_d)$ can be expressed in terms of the ``running"  coordinates $(y_1,\dots,y_d)$ by equations
$y_j-x_j=k_jx_j(y_d-x_d)$, $j=1,\dots,d-1$.
The square of the distance between  $(y_1,\dots, y_{d-1})$ and $M\Gamma$ is
$$
\sum\limits_{j=1}^{d-1}y_j^2=\sum\limits_{j=1}^{d-1}(x_j-k_jx_j(y_d-x_d))^2.
$$
The ``ordinate" $y_d={\mathcal C}\mu$ of the metacenter gives the minimum of the above expression and annihilates its derivatives (at $x_d=0$). Hence,
$$
\sum\limits_{j=1}^{d-1}k_jx_j(x_j-k_jx_jy_d)=0,\qquad\textrm{i.e.,}\qquad {\mathcal C}\mu=\frac{\sum\limits_{j=1}^{d-1}k_jx_j^2}{\sum\limits_{j=1}^{d-1}k_j^2x_j^2}\,.
$$
If $MT$ is the unit tangent vector to $\gamma$ at $M$, then, identifying $e_d^{\perp}$ with ${\mathbb R^{d-1}}$, writing $MT$   in spherical coordinates $\zeta=(\zeta_1,\dots,\zeta_{d-1})\in S^{d-2}$ and putting $(\zeta_1,\dots,\zeta_{d-1})=\frac{(x_1,\dots,x_{d-1})}{\sqrt{x_1^2+\dots+x_{d-1}^2}}$, we obtain (\ref{DVPkr1}).

{\bf Acknowledment}. The author is very thankful  to Mariangel Alfonseca, Alexander Fish, Carsten Sch\"utt,  Elisabeth Werner, Vlad Yaskin  and Ning Zhang for  very useful discussions. He is also very grateful  to Daniel Hug and Christos 
Saroglou for explaining several results and providing   references, and to Peter V\'arkonyi for pointing out the possible necessity of an extra condition in Lemma \ref{tr} for non-symmetric convex bodies.


\begin{thebibliography}{FSWZ}




\bibitem[A]{A}
{\sc H. Auerbach}, {\em Sur un probl\'eme de M. Ulam concernant l'\'equilibre des corps 
flottants}, Studia Mathematica
{\bf 7} (1938), no. 1, 121-142.

\bibitem[Al]{Al} {\sc Aleksandrov}, {\em Almost everywhere existence of the second differential of a convex function and some properties of convex surfaces connected with it} (in Russian), Uchenye Zapiski Len. Gos. Univ. Math. Ser. {\bf 6} (1939), 3-35.





\bibitem[B]{B} {\sc I. B\'ar\'any}, {\em Random  polytopes  in  smooth  convex  bodies},  Mathematika {\bf 39} (1992),  89-92; Corrigendum, Mathematika {\bf 51} (2004), 31.





\bibitem[BL]{BL}
{\sc  I. B\'ar\'any  and D. G. Larman},  {\em Convex bodies, economic cap coverings, random polytopes}, Mathematika, {\bf 35}(2) (1988), 274-291.


\bibitem[BF]{BF} {\sc H. Busemann and W. Feller}, {\em Krummungseigenschaften konvexer Fl\'ahen}, Acta. Math., {\bf 66}
(1935), 1-47.


\bibitem[BLW]{BLW}{\sc F. Besau, M. Ludwig and E. Werner}, {\em Weighted floating bodies and polytopal approximation},
Trans. Am. Math. Soc.,  {\bf 370}, (2018), 7129-7148 .


\bibitem[BMO]{BMO} {\sc  J.  Bracho, L. Montejano, D. Oliveros},  {\em Carousels,  Zindler curves and the floating body problem}, Per. Mat. Hungarica, vol. 49 {\bf 2} (2004), 9-23.

\bibitem[BP]{BP} {\sc H. Busemann and C. Petty}, {\em Problems on convex bodies}, Math. Scand. {\bf 4} (1956), 88--94. 



 \bibitem[CFG]{CFG}
{\sc H. T. Croft, K. J. Falconer, and R. K. Guy}, {\em Unsolved problems in geometry}, Problem Books in Mathematics,
Springer-Verlag, New York, 1991, Unsolved Problems in Intuitive Mathematics, II.


\bibitem[Da]{Da} {\sc A. Davidov}, {\em Theory of bodies floating in liquids in equilibrium},
Matser Thesis, Moscow State University, 1848, in 
{\em Life and works of A. Yu. Davidov},  (in Russian), by N.E. Zhukovsky, P.A. Nekrasov and P.M. Pokrovsky, 
Math. Sb. 1890, vol. 15, {\bf 1}, 1-57.

 \bibitem[DC]{DC} {\sc M. P. Do Carmo}, {\em Differential geometry of curves and surfaces}, IMPA, Rio de Janeiro, Brazil, ISBN 0-13-212589-7.



\bibitem[D]{D} {\sc C. Dupin}, {\em  Application de g\'eometie et de m\'echanique \'a la marine, aux ponts et chausse\'ees}, Paris, 1822.


\bibitem[DVP]{DVP} {\sc CH. J.  De La Vall\'ee Poussin}, {\em Lecons De  M\'ecanique Analytique}, Vol II, Paris, 1925 (in French), see also the Russian translation,  Moscow, 1949.



\bibitem[Fa]{Fa} {\sc K. J. Falconer}, {\em  Applications of a Result on Spherical Integration to the Theory of Convex Sets}, Amer. Math. Monthly, {\bf 90} (1983), 690-693.

\bibitem[FSWZ]{FSWZ} {\sc D. I. Florentin, K. Sch\"utt, E. M. Werner and N. Zhang}, {\em Convex floating bodies of Equilibrium}, arXiv:2010.09006.



\bibitem[Ga]{Ga} {\sc R. J. Gardner}, {\em Geometric tomography}, second ed., Encyclopedia of Mathematics and its Applications,
{\bf 58}, Cambridge University Press, Cambridge, 2006.

\bibitem[G]{G} {\sc E. N. Gilbert}, {\em How things 
float}, The American Mathematical Monthly {\bf 98} (1991), no. 3, 201-216.




\bibitem[Gr]{Gr} {\sc H. Groemer}, {\em Eine  kennzeichnende Eigenschaft der Kugel}, Eiseign. Math. (2) {\bf 7} (1961), 275-276.


\bibitem[H]{H} {\sc R. Howard}, {\em Alexandrov's theorem on the second derivatives of convex functions via Rademacher's theorem on the first derivatives of Lipschitz functions},
http://ralphhoward.github.io/SemNotes/Notes/alex.pdf

 \bibitem[HSW]{HSW} {\sc H. Huang, B. Slomka and E. Werner}, {\em Ulam floating bodies},  J. of London Math. Soc.,  {\bf 100} (2019), no. 2, 425-446.

 \bibitem[KO]{KO} {\sc A. Kurusa and T. \'Odor}, {\em Spherical floating bodies}, Acta Sci. Math. (Szeged), (2015), 81:3-4, 699-714.

\bibitem[L]{L} {\sc K. Leichtweiss}, {\em Zur Affinoberfl\"ache konvexer K\"orper}, Manuscripta Math., {\bf 56} (4), 429-464.

\bibitem[M]{M} {\sc R. D. Mauldin}, {\em The Scottish book, Mathematics from the Scottish Caf\'e with selected problems from the new Scottish book}, Second Edition, Birkh\"auser, 2015, ISBN 978-3-319-22896-9.


\bibitem[MP]{MP} {\sc V. D. Milman and A. Pajor}, {\em Isotropic position and inertia
ellipsoids and zonoids of the unit ball of a normed $n$-dimensional
space}, GAFA, Lecture
Notes in Math., {\bf 1376} (1989), Springer Berlin, 64-104.





\bibitem[Mo]{Mo} {\sc L. Montejano}, {\em On a  problem of Ulam concerning a characterization of the sphere}, Studies Appl. Math. {\bf 53} (1974), 243-248.

\bibitem[MR]{MR} {\sc M. Meyer and S. Reisner}, {\em   A geometric property of the boundary of symmetric convex bodies and convexity of flotation surfaces}, Geom. Dedicata {\bf 37} (1991), no. 3, 327-337.

\bibitem[MRS]{MRS} {\sc } {\sc S. Myroshnychenko, D. Ryabogin and C. Saroglou}, 
{\em  Star bodies with completely symmetric sections}, Int. Math. Res. Not., {\bf 10} (2019), 3015-3031. 

\bibitem[Na]{Na} {\sc F. Nazarov}, {\em Personal communication}, 
2010.

\bibitem[NSW]{NSW} {\sc S. Nagy, C. Sch\"utt and E. Werner}, {\em Data depth and  floating body}, Statistics Surveys {\bf 13} (2019), 52-118. 

\bibitem[Od]{Od} {\sc K. Odani}, {\em Ulam's floating body problem of two dimension}, Bull. of Archi Univ. of Education, {\bf 58} (2009), 1-4.

\bibitem[O]{O} {\sc S. P. Olovjanischnikoff}, {\em Ueber eine kennzeichnende Eigenschaft des Ellipsoides}, Leningrad State Univ. Ann. (Uchen. Zap.) {\bf 83} (1941), 113-128.

\bibitem[R1]{R1}{\sc D. Ryabogin}, {\em On a equichordal property for a pair of convex bodies}, arXiv: mathematics
2010.09864 

\bibitem[R2]{R2}{\sc D. Ryabogin}, {\em A negative answer to Ulam's Problem 19 from the Scottish Book}, arXiv: mathematics
1201.0393


\bibitem[Sch1]{Sch1}{\sc R. Schneider}, {\em Functional equations connected with rotations and their geometric applications}. L'Enseign. Math. {\bf 16} (1970), 297-305.


\bibitem[Sch2]{Sch2}{\sc R. Schneider}, {\em Convex Bodies: The Brunn-Minkowski theory}, Encyclopedia of
Mathematics and its Applications, Second expanded edition, 44, Cambridge University Press, Cambridge, 2014.


\bibitem[S1]{S1} {\sc C. Sch\"utt}, {\em On the affine surface area}, Proc. AMS., {\bf 118} (1993), 1213-18.

\bibitem[S2]{S2}{\sc C. Sch\"utt}, {\em Random polytopes and affine surface area}, Math. Nachr., {\bf 170} (1994), 227-249.




\bibitem[St]{St} {\sc A. Stancu}, {\em The floating body problem}, Bull. London Math. Soc., {\bf 38} (2006) 839–846.



\bibitem[SW1]{SW1} {\sc C. Sch\"utt and E. Werner}, {\em The convex floating body}, Math. Scand. {\bf 66} (1990), 275-290.

\bibitem[SW2]{SW2} {\sc C. Sch\"utt and E. Werner}, {\em Homothetic  floating body}, Geom. Dedicata. {\bf 49} (1994), 335-348.








\bibitem[T]{T} {\sc J. A. Thorpe}, {\em Elementary topics in Differential Geometry}, Underg. Texts in Math., Springer,  1979, ISBN 3-540-90357-7.


\bibitem[Tu]{Tu} {\sc E.  C. Tupper}, {\em An in Introduction to Naval Architecture} (Fifth Edition), 2013,  ISBN: 9780080982373.


\bibitem[U]{U} {\sc S. M. Ulam},  {\em A Collection of Mathematical Problems},  Interscience, New York, 1960, p. 38.



\bibitem[V1]{V1} {\sc P. L. V\'arkonyi},  {\em Floating body problems in two dimensions},  Stud. Appl. Math. {\bf 122} (2009), no. 2, 195–218.





\bibitem[V2]{V2} {\sc P. L. V\'arkonyi},  {\em Neutrally 
	floating objects of density $\frac{1}{2}$ in three dimensions}, Stud. Appl. Math. {\bf 130} (2013),
no. 3, 295-315.





\bibitem[Weg1]{Weg1} {\sc F. Wegner}, {\em Floating bodies of equilibrium},  Stud. Appl. Math. {\bf 111} (2003), no. 2, 167–183.

\bibitem[Weg2]{Weg2} {\sc F. Wegner}, {\em Floating bodies in  equilibrium in $2D$, the tire track problem and electrons in a parabolic magnetic fields}, arXiv:physics/0701241v3 (2007).





\bibitem[W]{W} {\sc E. Werner}, {\em Illumination bodies and affine surface area}, Stud. Math.,  {\bf 110}, (1994), 257-269.


\bibitem[Zh]{Zh} {\sc N. E. Zhukovsky}, {\em Classical mechanics},  Moscow, 1936 (in Russian).


\bibitem[Zi]{Zi} {\sc K. Zindler}, {\em \"Uber konvexe Gebilde II}, Monatsh. Math. Phys. {\bf 31} (1921), 25-57.

\end{thebibliography}
\end{document}